\documentclass[12pt]{article}
\usepackage{amssymb,amsmath,amscd}

\newtheorem{proposition}{Proposition}

\newtheorem{lemma}[proposition]{Lemma}
\newtheorem{remark}[proposition]{Remark}

\newtheorem{theorem}[proposition]{Theorem}

\def\+{{+\!\!\!+}}

\def\d{\partial}

\def\G{\Gamma}

\def\pmb#1{\setbox0=\hbox{#1}%
\kern.0em\copy0\kern-\wd0
\kern-.04em\copy0\kern-\wd0
\kern.08em\copy0\kern-\wd0
\kern-.04em\raise.0433em\box0 }         

\def\id{{\rm id}}
\DeclareMathOperator{\bis}{Bis}

\newcommand{\nc}{\newcommand}
\nc{\ber}{\begin{eqnarray}}
\nc{\eer}[1]{\label{#1}\end{eqnarray}}
\nc{\pek}[1]{\cite{#1}} \nc{\enr}[1]{(\ref{#1})}
\nc{\kal}[1]{{\cal{#1}}} \nc{\dott}{\;\cdot\;}
\nc{\Scl}{S_{\mathrm{cl}}}
\nc{\Diff}{{\mathit{Diff}}}

\def\G{{\cal G}}
\def\M{{\cal M}}

\def\sphere{{\mathbb S}}
\def\Z{{\mathbb Z}}
\def\R{{\mathbb R}}
\def\C{{\cal C}}

\newcommand{\be}{\begin{equation}}
\newcommand{\ee}{\end{equation}}
\newcommand{\bea}{\begin{eqnarray}}
\newcommand{\eea}{\end{eqnarray}}

\newcommand{\Section}[1]{\section{#1} \setcounter{equation}{0}}
\begin{document}
\begin{center}

\vskip .3in \noindent

\vskip .1in

{\Large \bf{Geometric quantization and non-perturbative Poisson sigma model}}
\vskip .2in

{\bf Francesco Bonechi}$^a$\footnote{e-mail address: Francesco.Bonechi@fi.infn.it},
{\bf Alberto S. Cattaneo}$^b$\footnote{e-mail address: alberto.cattaneo@math.unizh.ch}
 and  {\bf Maxim Zabzine}$^{c}$\footnote{e-mail address: m.zabzine@qmul.ac.uk} \\


\vskip .15in

\vskip .15in
$^a${\em I.N.F.N. and Dipartimento di Fisica}\\
{\em  Via G. Sansone 1, 50019 Sesto Fiorentino - Firenze, Italy} \\
\vskip .15in
$^b${\em   Institut f\"ur Mathematik, Universit\"at Z\"urich-Irchel}\\
{\em Winterthurerstrasse 190, CH-8057 Z\"urich, Switzerland}
\vskip .15in
$^c${\em School of Mathematical Sciences, Queen Mary, University of London}\\
{\em Mile End Road, London E1 4NS, UK}

\bigskip


 \vskip .1in
\end{center}
\vskip .4in

\begin{center} {\bf ABSTRACT }
\end{center}
\begin{quotation}\noindent
  In this note we point out the striking relation between the conditions
 arising within geometric quantization and the non-perturbative Poisson
 sigma model. Starting from the Poisson sigma model, we analyze necessary requirements
 on the path integral measure which imply a certain
 integrality condition for the  Poisson cohomology class $[\alpha]$.
 The same condition was considered before by Crainic and Zhu but in a different context.
 In the  case when $[\alpha]$ is in the image of the sharp map
  we reproduce the Vaisman's condition for prequantizable Poisson manifolds.
  For integrable Poisson manifolds we show, with a different procedure than in Crainic and Zhu,
  that our integrality condition implies the prequantizability of the
  symplectic groupoid. Using the relation between prequantization and symplectic reduction
   we construct the explicit prequantum line bundle for a symplectic  groupoid.
   This picture supports the program of quantization of Poisson manifold via symplectic
    groupoid.
 At the end we  discuss the case of a generic coisotropic D-brane.
\end{quotation}
\vfill
\eject


\section{Introduction}
\label{s:intro}

 Quantization is generally understood as a transition from classical to quantum
 mechanics. In mathematics  a quantization of Poisson manifold should be  a  prescription
  able to produce a structure which physicists would agree to call the quantum theory
 associated with the classical system given by the Poisson manifold.
 However this prescription is far from being unique and
 different people give a different meanings to the word
{\it quantization}. In this note we are interested in two different incarnations
 of the word {\it quantization}: {\it geometric quantization}\/ and {\it deformation quantization}.

 {\it Deformation quantization}\/ deals with formal associative deformations of $C^\infty(M)$.
In \cite{Kontsevich:1997vb} Kontsevich gave a general formula for the {\it deformation quantization}\/
  of the algebra of functions on a Poisson manifold $(M, \alpha)$.
 Later in \cite{Cattaneo:1999fm} it was established that the perturbative path integral expansion
 of the Poisson sigma model over
 the two-dimensional disk $D$ leads to the Kontsevich's star product.

The Poisson sigma model, introduced in \cite{Ikeda:1993fh, Schaller:1994es}, is a topological
 two-dimensional field theory  with target a Poisson manifold $M$, whose
  Poisson tensor we will denote by $\alpha$ throughout. Let $\Sigma$ be a two-dimensional oriented compact manifold
 with a boundary. The starting point is the classical action functional $S$
 defined on the space of vector bundle morphisms $\hat{X}\colon  T\Sigma \rightarrow T^*M$
 from the tangent bundle $T\Sigma$ to the cotangent
 bundle $T^*M$ of the Poisson manifold $M$.
 Such a map $\hat{X}$ is given by its base map $X\colon  \Sigma \rightarrow M$  and the linear map $\eta$ between
 fibers, which may also be regarded as
 a section in $\Gamma(\Sigma, Hom(T\Sigma, X^*(T^*M)))$.
  The pairing $\langle\,\,,\,\,\rangle$ between the cotangent and tangent space at each point of $M$
  induces a pairing between the differential forms on $\Sigma$ with values in the pull-backs
 $X^*(T^*M)$ and $X^*(TM)$ respectively. It is defined
 as pairing of the values and the exterior product of differential forms.
 Then the action functional $S$ of the theory is
\begin{equation} S(X,\eta) = \int\limits_\Sigma  \langle \eta,  dX\rangle +
\frac{1}{2} \langle \eta, (\alpha \circ X) \eta \rangle .
\label{definPS}
\end{equation}

 Here $\eta$ and $dX$ are viewed as one-forms on $\Sigma$ with the values in the pull-back of
 the cotangent and tangent bundles of $M$ correspondingly.
 Thus, in local coordinates, we can rewrite the action (\ref{definPS}) as follows:
 \begin{equation}
 S(X,\eta) = \int\limits_D  \eta_\mu \wedge dX^\mu + \frac{1}{2} \alpha^{\mu\nu}(X) \eta_\mu
 \wedge \eta_\nu .
\label{definPSloccoor}\end{equation}
 The variation of the action gives rise to the following equations of
 motion
\begin{equation}
 d\eta_\rho + \frac{1}{2} (\d_\rho \alpha^{\mu\nu}) \eta_\mu \wedge \eta_\nu =0,\,\,\,\,\,\,\,\,\,\,\,
 dX^\mu + \alpha^{\mu\nu}\eta_\nu = 0 .
\label{eqmotion}\end{equation}
 In covariant language these equations are equivalent to the statement that
 the bundle morphism $\hat{X}$ is  a Lie algebroid morphism  from
 $T\Sigma$ (with standard Lie algebroid structure) to $T^*M$ (with Lie algebroid structure canonically
 induced by the Poisson structure).
The action (\ref{definPSloccoor}) is invariant under the infinitesimal gauge transformations
\begin{equation}
 \delta_\beta X^\mu = \alpha^{\mu\nu} \beta_\nu,\,\,\,\,\,\,\,\,\,\,\,\,\,
 \delta_\beta \eta_\mu = - d\beta_{\mu} - (\d_\mu \alpha^{\nu\rho}) \eta_\nu \beta_\rho ,
\label{gaugetransf}\end{equation}
 which form a closed algebra only on-shell (i.e., modulo the equations of motion (\ref{eqmotion})).
 We are interested in the situation when $\d \Sigma \neq \emptyset$.
 Following \cite{Cattaneo:1999fm} we first consider the boundary conditions
\begin{equation}
\eta_{t}|_{\d \Sigma}=0\,\,\,\,\,\,\,\,\,\,\,\,\,\,
\beta|_{\d \Sigma}=0 \label{CFboundary}\end{equation}
 where  $t$ corresponds to the direction tangent to the boundary.
 More general boundary conditions will be
 discussed in Section \ref{s:boundary}.

In \cite{Cattaneo:1999fm} it has been shown that the star product is given by the semiclassical
 expansion of the path integral of the Poisson sigma model over a disk $D$
\begin{equation}
 f * g (x) = \int\limits_{X(\infty)=x} f(X(1)) g(X(0)) e^{\frac{i}{\hbar} S(X,\eta)} dX d\eta ,
\label{defpath}\end{equation}
 where $0,1, \infty$ are any three cyclically ordered points on the unit circle $\d D$.
 The semiclassical expansion is to be understood as as an expansion around the trivial
 classical solution $X=x$ and $\eta=0$.  As it stands the integral (\ref{defpath}) is not
  well-defined due to the gauge symmetries and renormalization. However
 at perturbative level this can be fixed \cite{Cattaneo:1999fm}.

In this note we explore the idea that the formula (\ref{defpath}) could make sense
 also outside of the perturbative expansion around the trivial solution.
 Although we do not know how to define
 the path integral non-perturbatively in general, we can try to do different
  consistency checks.
 Thus, proposing different consistency tests for the Poisson sigma model on the sphere or on the disk
 with boundary conditions \eqref{CFboundary}, we arrive to
 the following integrality condition, which is the main result of the paper:
\begin{equation}
\frac{1}{2\pi \hbar} S (X,\eta)
\,\,\in\,\,\Z,
\label{stateresultS}\end{equation}
for every classical solution $X,\eta$ on the sphere or on the disk.
We also show that this condition is equivalent to having,
for every two-cycle $c_2$
 which is the image of the base map of
 a Lie algebroid  morphism $T \sphere^2 \rightarrow T^*M$,
\begin{equation}
 \frac{1}{2\pi\hbar} \int\limits_{c_2}\,\, \omega_{\cal L}\,\,\in\,\,\Z
\label{stateresult}\end{equation}
 where $\omega_{\cal L}$ denotes the induced symplectic form
on the symplectic leaf ${\cal L}$ that contains $c_2$.
 This integrality condition is related to
the different integrality conditions which appear within the {\it geometric quantization
 program}. It is also the necessary and sufficient condition for the symplectic groupoid
 of the given Poisson manifold to be prequantizable \cite{crainic1}, a result that we will rederive,
 using different methods, in Section~\ref{s:groupoid}.
 The principal aim of this paper is to give different derivations
 of condition (\ref{stateresult}) and to relate it to the known integrality conditions.
 It appears that our considerations give strong support to the program \cite{Weinstein}
 of quantizing Poisson manifolds via their corresponding symplectic groupoids.
It seems plausible that considering the Poisson sigma model on any two-manifold $\Sigma$ should yield
a stronger integrality condition; viz., \eqref{stateresult} should be satisfied for any $c_2$ which is
 the image of the base map of a Lie algebroid
 morphism $T \Sigma \rightarrow T^*M$ for any $\Sigma$. However, we still do not have a proof of this fact.\footnote{
 The first attempt to discuss the role of integrality in the context
 of PSM has been done in \cite{Schaller:1995xk}.}
 
Moreover this integrality condition has a cohomological meaning. On a classical 
solution $\hat{X}=(X,\eta):T\Sigma\rightarrow T^*M$ ($\partial\Sigma=\emptyset$) the action functional can be written as
\begin{equation}
 \Scl  (X,\eta)= - \frac{1}{2} \int\limits_{\Sigma} \langle \eta, (\alpha \circ X) \eta\rangle \;.
\label{actonsol}\end{equation}
 This expression depends only on the cohomology class of $[\alpha] \in
  H^\bullet_{LP}(M,\alpha)$ and the homology class of $(\Sigma, \hat{X})$ as explained in details in the 
  Appendix; in fact (\ref{actonsol}) can be reintepreted in terms of a natural paring $\ll\,\,,\,\,\gg$ 
  between the Poisson cohomology $H^\bullet_{LP}(M,\alpha)$ and the homology whose cycles are algebroid
   morphisms. In particular it makes sense to talk of integral Poisson 
  tensor, which corresponds to requiring that
    \begin{equation}
     \ll \alpha, (c_2, \hat{X})
   \gg = 2\pi n \hbar,\,\,\,\,\,\,\,\,\,\,\,n \in \mathbb{Z},
    \label{newcond123}\end{equation}
 for any algebroid morphism $\hat{X}:Tc_2\rightarrow T^*M$, and, more generally, of integral Poisson 
 cohomology, so that
    \begin{equation}
  \left[\frac{1}{2\pi \hbar}\alpha\right] \in H^2_{LP}(M, \alpha, \mathbb{Z}).
  \label{inegerPoisco}\end{equation}


The organization of our presentation is as follows: In Section \ref{s:geometric}
  we review the relevant concepts
 from the geometric quantization. In particular we discuss the known integrability conditions on
  symplectic forms and Poisson tensors arising in geometric quantization.
 In  Section \ref{s:gauge} we consider the on-shell gauge transformations and the value of the
 action functional on classical solutions of the Poisson sigma model. We derive a generalized
 integrality condition for a Poisson tensor in subsection \ref{subs:onshell}.
  In subsection \ref{subs:symplectic} we interpret this condition in the symplectic case
 and discuss some peculiarities of the Poisson sigma model over symplectic manifolds.
 In subsection \ref{subs:poisson} we return to the case of general Poisson manifolds,
 discuss the cohomological nature of our integrality condition and recover Vaisman's condition
 for the prequantizability of the Poisson manifold as a special case.
 In Section \ref{s:groupoid} the integrality condition is rederived in a different context:
 Namely, we study the relation between symplectic reduction and prequantization
 in the context of the Poisson sigma model and its reduced phase space.
 Thus, assuming the integrality condition  for an integrable Poisson manifold, we
 show that the corresponding symplectic groupoid is prequantizable.
 In Section \ref{s:boundary} general boundary conditions are considered and
 the corresponding integrality conditions are found.
 In Section \ref{s:end} we present concluding remarks with a possible interpretation
 of our results. In  the Appendix we collect the technicalities regarding the relevant (co)homology
  groups.

\section{Geometric quantization}
\label{s:geometric}

This Section is just a reminder of some relevant elements of geometric
 quantization. Namely, we are going to review the  necessary and
 sufficient conditions for the existence of prequantization bundles for symplectic and Poisson
 manifolds.

 The geometric quantization theory of Kostant \cite{Kostant} and Souriau \cite{Souriau} was first developed
 for symplectic manifolds and then further generalized to Poisson manifolds.
 Within this approach the quantization is done in two steps: prequantization
 which gives a linear representation of $(C^\infty(M), \{\,\,,\,\,\})$ by operators
 on a complex vector space and then quantization where one restricts to a convenient
 subalgebra of $(C^\infty(M), \{\,\,,\,\,\})$ (for review see \cite{Kirillov, Vaisman, Woodhouse}).

 We are interested in the first step of this construction. Modulo a certain obstruction,
 the prequantization problem can be solved by considering the space of sections of a complex
 line bundle $L \rightarrow M$, the prequantization bundle and a well-chosen prequantization
 formula. In other words, to each $f \in C^{\infty}(M)$ there corresponds an operator $\hat{f}$
 acting on $\Gamma(L)$ such that the map $f \mapsto \hat{f}$ is linear and
\begin{equation}
 - i \hbar\widehat{ \{f, g \}} =  [\hat{f}, \hat{g}] \equiv  \hat{f} \hat{g} - \hat{g}\hat{f}
\label{definprincpsl}\end{equation}
 with $1 \rightarrow \hat{1}$.

 Let us start from the symplectic case. Assuming that the line bundle $L$ exists, we can
 define $\hat{f}$ by
\begin{equation}
 \hat{f} s =  - i\hbar \nabla_{X_f} s + f s,
\label{symplectcase}\end{equation}
 where $\nabla$ is a covariant derivative on $L$ and $X_f$  is the Hamiltonian vector
 field of $f$ (i.e., $X_f g = \{f, g\}$ for every function $g$).
  Then, (\ref{definprincpsl}) is equivalent to
\begin{equation}
 c_1 (L) = - \frac{1}{2\pi\hbar} [\omega] ,
\label{firstCh}\end{equation}
 where $c_1(L)$ is the first Chern class of the line bundle $L$. Since the first Chern class
 is integral (i.e., $c_1(L) \in H^2(M, \mathbb{Z})$), one arrives at the following quantization condition
 for the symplectic form:
\begin{equation} [\frac{1}{2\pi \hbar} \omega] \in H^2(M, {\mathbb Z}).
\label{quantizat}\end{equation}
 The converse is also true. Namely, if the symplectic form satisfies (\ref{quantizat})
 then there exists the line bundle $L$ with connection $\nabla$ such that (\ref{definprincpsl})
 and (\ref{symplectcase}) are fulfilled.
 Observe that the representation $f\mapsto\hat f$ is faithful.

In the Poisson case one can generalize the above construction \cite{Vaismanpap}.
The idea is that it is enough to define $\nabla$ only along the symplectic leaves
 and thus one needs only the partial connection on $L$. The right concept is
 the contravariant derivative for the vector bundle $L$ over a Poisson
 manifold $(M, \alpha)$:
\begin{equation}
 D_w \colon \,\,\Gamma(L)\,\,\rightarrow\,\,\Gamma(L),\,\,\,\,\,\,\,\,\,w\in T^*M,
\label{defcont}\end{equation}
 such that for any $f\in C^{\infty}(M)$ and $s\in \Gamma(L)$
\begin{equation} D_{fw} s = f D_w S,\,\,\,\,\,\,\,\,\,\,\,\,\,\,\,D_w fs =f
D_w s + w (\sigma f) s, \label{defcontder}\end{equation}
 where $\sigma = [\alpha,\,\,\,\,]_s$
 is the Schouten--Nijenhuis bracket with $\alpha$.

 The sharp map  $\sharp\colon  T^*M \rightarrow TM$ is defined by $w(\sharp p) =
 \alpha(p, w)$ with $w,p \in T^*M$.
  It induces a homomorphism on the relevant cohomologies:
\begin{equation}
 \sharp\colon \,\,\,\,\,H_{deR}(M, \mathbb{R})\,\,\rightarrow\,\,H_{LP}(M,\alpha) ,
\label{definsharpco}\end{equation}
 where $H_{deR}(M, \mathbb{R})$ is de~Rahm cohomology on forms and
 $H_{LP}(M,\alpha)$ is Lichnerowicz--Poisson cohomology on contravariant
 antisymmetric tensors with the differential $\sigma$.

Using a contravariant derivative the prequantization formula becomes
\begin{equation}
 \hat{f}s = - i\hbar D_{df} s + fs
\label{poipre}\end{equation}
 and condition (\ref{definprincpsl}) implies
 that
\begin{equation}
 C_D = - \frac{i}{\hbar} \alpha,
\label{curvature}\end{equation}
 where $C_D$ is the curvature of D,
\begin{equation} C(w,p)s = D_w D_p s - D_p D_w s + D_{[w,p]} s,
\label{defcurvature}\end{equation}
 with $[\,\,,\,\,]$ the Koszul brackets on one-forms of a Poisson manifold.
 Thus, the Poisson-Chern class of $L$ is
\begin{equation}
 pc_1(L) = - \frac{1}{2\pi \hbar}[\alpha] ,
\label{PoisCh}\end{equation}
 which is the image of the real Chern class $c_1(L)$ under (\ref{definsharpco}).
 So
  the existence of a prequantization bundle requires that
 the preimage of $\frac{1}{2\pi \hbar}[\alpha]$ under (\ref{definsharpco}) should be an element of
 $H^2(M, \mathbb{Z})$ (i.e., $H^2(M,\mathbb{Z})$ is sent by
 the inclusion $\mathbb{Z} \subset\mathbb{R}$ to $H^2_{deR}(M, \mathbb{R})$).
 The converse is also true. Thus, this construction gives a representation of
 $C^\infty(M)$, which is however not always faithful.
 For further details and relevant concepts, see textbook by Vaisman \cite{Vaisman}.

 To summarize, a Poisson manifold $(M, \alpha)$ has a prequantization bundle if{f}
  there exist a vector field $v$  and a closed two-form $w$ that
  represents an integral cohomology class of $M$ (i.e., $\frac{1}{2\pi \hbar} [w]\in  H^2(M, \mathbb{Z})$),
  such that
\begin{equation}
 \alpha = - \sharp w + {\cal L}_v \alpha.
\label{defpeqabnd}\end{equation}
 Within the geometric quantization framework such $(M, \alpha)$ is called
 quantizable. 
   
We finally show that if the Poisson tensor satisfies (\ref{defpeqabnd}) then it also satisfies (\ref{stateresultS}). Let $(X,\eta)$ be a solution of (\ref{definPSloccoor}), for a generic surface $\Sigma$. By direct computation we get
the value of the action (\ref{definPS}) on $(X,\eta)$ as
\begin{equation}
 \Scl (X, \eta) = \int\limits_{\Sigma} X^*(w) + \int_{\partial\Sigma} i_v\eta 
\label{solvalue}\end{equation}
 for the generic boundary condition. If we specialize to $\Sigma = \sphere^2$, there is no boundary term and
the integrality of $w$ obviously implies (\ref{stateresultS}). However, the opposite is obviously not true, since condition (\ref{stateresultS}) makes sense also when $[\alpha]$ is not in the image of the sharp map.

\section{Gauge transformations of the Poisson sigma model}
\label{s:gauge}

In this Section we analyze the role of the integrality condition
(\ref{stateresult}) for the Poisson tensor $\alpha$ in the
nonperturbative definition of the Kontsevich formula
(\ref{defpath}).

There is a simple argument supporting the need of
(\ref{stateresult}). The Kontsevich formula appears from the
perturbative expansion around the classical trivial solution
$X=x,\eta=0$ of (\ref{defpath}). One can ask if there exist
instantons, e.g. inequivalent nontrivial solutions, around which
to expand the integral. The expansion around a classical solution
$(X,\eta)$ has the form of an asymptotic series
\begin{equation}
 e^{\frac{i}{\hbar}\Scl (X,\eta)} (b_0 + \hbar b_1 + \hbar^2 b_2 + ...),
\label{expsol}
\end{equation}
and in general is not just a power series. However, due to
boundary conditions (\ref{CFboundary}) every solution maps the
boundary $\d D$ to a single point $x$ and defines an algebroid
morphism from $T\sphere^2$ to $T^*M$. Therefore property
(\ref{stateresultS}) implies that the action evaluated on a solution is
equal to $2\pi n \hbar$ and there is no exponential factor in the
expansion (\ref{expsol}). This ensures that the full
expansion (i.e., over all non-trivial solutions) of (\ref{defpath})
will give rise anyway to a formal power series in $\hbar$ (which is
correct from the point of view of star
 products).

In this section we investigate  the role of the integrality
condition in the nonperturbative definition of (\ref{defpath}). In
subsection \ref{subs:onshell} we analyze the on-shell gauge
transformations of the PSM for a generic closed surface $\Sigma$
and for the disk D with boundary conditions (\ref{CFboundary}). In
the expansion of the path integral (\ref{defpath}) around a
non-trivial solution, we have to insure that the exponential of
the classical action $\Scl$ (i.e., the value of the action
functional (\ref{definPS}) on a classical solution)
 is well-defined. In other words, if two Lie algebroid morphisms $\hat{X}_1 = (X_1,\eta_1)$
 and $\hat{X}_2 = (X_2, \eta_2)$ are gauge equivalent, then we must require the following condition:
\begin{equation}
 e^{\frac{i}{\hbar} \Scl (X_1,\eta_1)} = e^{\frac{i}{\hbar} \Scl (X_2,\eta_2)}.
\label{valclase}
\end{equation}
There are several ways of integrating the infinitesimal on-shell
transformations (\ref{gaugetransf}). The integrality condition
(\ref{stateresult}) comes as a necessary condition to satisfy
(\ref{valclase}) in the PSM on the disk with boundary conditions
(\ref{CFboundary}) if we choose the finite transformations to
define a groupoid action of $\G^{D,M}$---see \eqref{gdm} and
\eqref{gauge group} for its definition. In subsections
\ref{subs:symplectic} and \ref{subs:poisson} we perform a formal
analysis of the path integral formula (\ref{defpath}). The picture
that comes out is consistent with this choice of on--shell gauge
transformations; moreover, the integrality condition
(\ref{stateresult}) allows us to reduce (\ref{defpath}) to a
quantum mechanical path integral on the symplectic groupoid
integrating $M$. For pedagogical reasons, in subsection
\ref{subs:symplectic} we analyze first the symplectic case where
these ideas emerge in a simple way; the analysis is repeated in
subsection \ref{subs:poisson} for the general Poisson case, at the
cost of a more technical discussion involving the underlying
algebroid structures.

Finally observe that we have assumed the Poisson manifold to be
integrable. If this is not the case, a topological groupoid
$\G(M)$ ``integrating'' $M$ exists anyway as the reduced phase
space of the Poisson sigma model \cite{sympgroup}. This groupoid
may also be regarded as a stacky Lie groupoid \cite{TZ} and
possibly the considerations of the present Section go through.
This would imply the integrality condition \eqref{stateresult}
also in the nonintegrable case (notice that the condition is
well-defined for every Poisson manifold).

\subsection{On-shell gauge transformations}
\label{subs:onshell}

 We start by considering the Poisson sigma model defined over a closed surface $\Sigma$.

A Lie algebroid $A\rightarrow M$ is a vector bundle $A$ over $M$ with a Lie algebra
 structure (over $\mathbb R$)  on the sections of $A$ and a vector bundle morphism $\rho \colon  A
 \rightarrow TM$ (called the anchor) which induces a Lie algebra homomorphism from sections
  of $A$ to vector fields on $M$. The bracket on sections of $A$ is required to satisfy
   $[f e_1, e_2] = f [e_1, e_2] - (\rho(e_2) f )e_1$ for all $f \in C^\infty (M)$ and
   $e_1,e_2\in\Gamma(A)$.
 A Poisson manifold $(M,\alpha)$ induces a Lie algebroid structure on the cotangent bundle $T^*M$
 with $\sharp$ as its anchor; as for the Lie bracket on its sections (a.k.a.\ the Koszul bracket),
 it is enough to define it on
 exact one-forms by $[df, dg] \equiv d\{f, g\}$. The tangent bundle $T\Sigma$
 carries a canonical Lie algebroid structure with anchor the identity and the standard Lie
 bracket. The bundle morphism
\begin{equation}
\begin{array}{lll}
 T\Sigma & \stackrel{\eta}{\longrightarrow} & T^*M \\
\Big\downarrow && \Big\downarrow \\
 \Sigma & \stackrel{X}{\longrightarrow} & M
\end{array}
\label{lielagmorsj}
\end{equation}
satisfies the equations of
motion (\ref{eqmotion}) if{f} $(X,\eta)$ is a Lie algebroid morphism.
 In other words, the action functional, which is defined for any bundle morphism, is extremized
 when the morphism is a Lie algebroid morphism. This result follows straightforwardly
 from the definition of Lie algebroid morphisms given in \cite{higgins}.
 Thus, the space of classical solutions of the Poisson sigma model is the space ${\rm Mor}(T\Sigma,T^*M)$
  of Lie algebroid morphisms.

 A Lie groupoid ($\G \rightrightarrows M, s, t)$ is a smooth manifold $\G$ with surjective submersions $s$ and $t$
  from $\G$ to $M$ and a smooth multiplication map from $\G^{(2)}=\{
   (g, h) \in \G \times \G | t(h)=s(g)\}$ to $\G$ making $\G$ into  a category in which all elements
    are invertible. The points in $M$ are simultaneously considered as the objects of the
     category and the identity morphisms.  A Lie groupoid is called source-simply-connected
      (ssc) if the $s$-fibers are connected and simply connected.
          The vector bundle $\ker(d s)|_M$ has a natural structure of a Lie algebroid over $M$
      with anchor $dt$ and Lie bracket induced by multiplication. We denote by $A(\G)$
       the Lie algebroid of the Lie groupoid $\G$. Not all Lie algebroids arise in this way.
        Those which do are called integrable.

        A morphism of Lie groupoids induces a morphism of the corresponding Lie algebroids
         by taking its differential at the identity sections.  Let us  recall the second Lie
         theorem for Lie algebroids:

\begin{theorem}Let $\G$ and $\tilde{\G}$ be two Lie
groupoids integrating the Lie algebroids $A(\G)$ and
$A(\tilde{\G})$, respectively. If $\G$ is source simply connected
(ssc), then for any Lie algebroid morphism $A(\G)\rightarrow
A(\tilde{\G})$ there exists a unique integrating Lie groupoid morphism
$\G\rightarrow\tilde{\G}$.
\end{theorem}

For more information on Lie algebroids and
         groupoids, see \cite{makenzie, silva, moerdijk}.

A Poisson manifold is said to be integrable if the corresponding Lie algebroid $T^*M$ is integrable.
If this is case, we will denote the corresponding ssc Lie groupoid simply by ${\cal G}(M) $.
Indeed, ${\G}(M)$ carries a symplectic
  structure which is compatible with groupoid multiplication and hence it is called a symplectic
   groupoid.
    Let $\Pi(\Sigma) \rightrightarrows \Sigma$ be the
    fundamental groupoid of $\Sigma$ which is defined
    as follows: $\Pi(\Sigma)$ consists of homotopy classes $[c_{uv}]$ of curves $c_{uv}$
    in $\Sigma$ starting at $v\in \Sigma$ and ending at $u\in \Sigma$, with $s([c_{uv}])=
    v$ and $t([c_{uv}])=u$, and its multiplication law is induced by concatenation.
Observe that for $\Sigma$ simply connected
    $\Pi(\Sigma)=\Sigma\times\Sigma$.
    Since $\Pi(\Sigma)$ is the ssc groupoid integrating $T\Sigma$, we can associate to any
   Lie algebroid morphism (\ref{lielagmorsj}) a unique Lie groupoid morphism
\begin{equation}
\begin{array}{lll}
  \Pi(\Sigma) & \stackrel{\cal X}{\longrightarrow} &\G(M)\\
{\scriptstyle s}\Big\downarrow\Big\downarrow {\scriptstyle t}
 &&{\scriptstyle s}\Big\downarrow\Big\downarrow {\scriptstyle t} \\
  \,\,\,\,\Sigma & \stackrel{X}{\longrightarrow} &\,\,\, M
\end{array}
\label{groupmorp}\end{equation}
 Therefore we can identify the space ${\rm Mor}(T\Sigma,T^*M)$ of Lie algebroid
morphisms with the space ${\rm Mor}(\Pi(\Sigma),\G(M))$ of Lie
groupoid morphisms.

For the case of integrable Poisson manifolds $(M,\alpha)$, we may then describe the space
 of classical solutions of the Poisson sigma model as ${\rm Mor}(\Pi(\Sigma),\G(M))$.
  Next we would like to introduce the gauge group which naturally acts on
   this space (i.e., which sends solutions to solutions).

 Let us consider the infinite-dimensional groupoid $\G^\Sigma =
\{\hat{\Phi}\colon \Sigma\rightarrow\G(M)\}$ over
$M^\Sigma=\{\Phi\colon \Sigma\rightarrow M\}$ with structure maps defined pointwise.
Namely, we define
source and target by  $s(\hat{\Phi})(u)=s({\hat\Phi}(u))$,
$t(\hat{\Phi})(u)=t({\hat\Phi}(u))$ for $u \in \Sigma$ and multiplication by
 $\hat{\Phi}_1\hat{\Phi}_2 (u) =
\hat\Phi_1(u)\hat\Phi_2(u)$.   A section $\lambda$ of the
associated algebroid $A(\G^\Sigma)$
(this algebroid has been defined, at least for $\Sigma$ one-dimensional, in \cite{BC} intrinsically in terms
of the Lie algebroid $T^*M$, so it exists also for nonintegrable Poisson manifolds)
is defined by giving
 a section $\lambda(\Phi) \in \Gamma(\Phi^*T^*M)$ for every $\Phi\in
M^\Sigma$. There is a natural groupoid action\footnote{This action
may actually be extended to the space of maps
$\Pi(\Sigma)\to\G(M)$ that are compatible with $s$ and $t$ maps
(i.e., maps for which \eqref{groupmorp} is a commutative diagram),
which however are not required to be compatible with
multiplication.} of $\G^\Sigma$ on ${\rm Mor}(\Pi(\Sigma),\G(M))$
 which is given by
\begin{equation}
X_{\hat\Phi}(u)=t(\hat{\Phi})(u)\,,\,\,\,\,\,\,\,\,\,\,\,\,\,\,\,\,\,
{\cal X}_{\hat{\Phi}}([c_{uv}])={\hat{\Phi}}(u){\cal
X}([c_{uv}])\hat{\Phi}(v)^{-1},
\label{gauge_trans_gpd}\end{equation} where $(X, {\cal X}),
(X_{\hat\Phi}, {\cal X}_{\hat\Phi})\in{\rm
Mor}(\Pi(\Sigma),\G(M))$,
 $\hat{\Phi}\in\G^\Sigma$ with $s(\hat{\Phi})=X$  and $[c_{uv}]$
 is the homotopy class of a curve $c_{uv}$ in $\Sigma$.
 From (\ref{gauge_trans_gpd}) one can easily deduce the action
 of $\G^\Sigma$ on the space of Lie algebroid morphisms ${\rm Mor}(T\Sigma,T^*M)$,
  $(X,\eta)\rightarrow (X_{\hat\Phi}, \eta_{\hat\Phi})$. However the
    concrete expression is not relevant for our discussion.

One could also consider a group action by introducing the group of
bisections of $\G^\Sigma$,
\begin{equation}
\bis(\G^\Sigma)=\{\sigma\colon M^\Sigma\rightarrow\G^\Sigma\,|\,
s\circ\sigma=\id,\ t\circ\sigma=\psi_{\sigma}\in
\Diff(M^\Sigma)\}. \label{defbis}\end{equation} It is clear that
formula (\ref{gauge_trans_gpd}) defines in a straightforward way a
group action of $\bis\G^\Sigma$ on ${\rm Mor}(\Pi(\Sigma),\G(M))$.
The orbits of $\G^\Sigma$ contain those of $\bis\G^\Sigma$ and
they coincide only when it is true that for every
$\hat\Phi\in\G^\Sigma$ there exists a bisection $\sigma$ passing
through it, e.g. $\sigma(s({\hat\Phi}))=\hat\Phi$. Moreover let us
remark that $(\bis\,\G)^\Sigma=\{\hat{\sigma}\colon
\Sigma\rightarrow\bis\,\G(M)\}$ is a subgroup of
$\bis(\G^\Sigma)$. Namely, there is a Lie group morphism
\[
\Psi\colon\begin{array}[t]{ccc}
 (\bis\,\G)^\Sigma&\rightarrow&\bis(\G^\Sigma)\\
\hat\sigma&\mapsto&\sigma
\end{array}
\]
with $\sigma(\Phi)(u)\equiv\hat\sigma(u)(\Phi(u))$, $u\in\Sigma$,
$\Phi\in M^\Sigma$. We also have a left-inverse map
$\widetilde\Psi\colon\sigma\mapsto\hat\sigma$ defined by
$\hat\sigma(u)(m)\equiv\sigma(\underline m)(u)$, $u\in\Sigma$,
$m\in M$, where $\underline m$ denotes the constant map with value
$m$. Thus, $\Psi$ is injective.

There is an interesting subset of gauge transformations which can
 be described as follows.
 For every map $f\colon\Sigma\rightarrow\Sigma$, its differential $f_*\colon  T\Sigma\rightarrow T\Sigma$
 is  a Lie algebroid morphism.
 We denote by
$ {\cal F} \colon \Pi(\Sigma)\rightarrow\Pi(\Sigma)$
the corresponding Lie groupoid morphism.
 Thus, using a map $f$, any
groupoid morphism $(X, {\cal X})$ can be transformed to a new groupoid morphism $(X\circ
f,{\cal X}\circ {\cal F})$ as follows:
\begin{equation}
\begin{array}{lll}
  \Pi(\Sigma) & \stackrel{{\cal X} \circ {\cal F}}{\longrightarrow} &\G(M)\\
{\scriptstyle s}\Big\downarrow\Big\downarrow {\scriptstyle t}
 &&{\scriptstyle s}\Big\downarrow\Big\downarrow {\scriptstyle t} \\
  \,\,\,\,\Sigma & \stackrel{X \circ f}{\longrightarrow} &\,\,\, M
\end{array}
 =
\begin{array}{lllll}
\Pi(\Sigma) & \stackrel{\cal F}{\longrightarrow}
    &\Pi(\Sigma) & \stackrel{\cal X}{\longrightarrow} &\G(M)\\
{\scriptstyle s}\Big\downarrow\Big\downarrow {\scriptstyle t}
&&{\scriptstyle s}\Big\downarrow\Big\downarrow {\scriptstyle t}
 &&{\scriptstyle s}\Big\downarrow\Big\downarrow {\scriptstyle t} \\
  \,\,\,\,\Sigma & \stackrel{f}{\longrightarrow}
 & \,\,\,\,\Sigma & \stackrel{X}{\longrightarrow} &\,\,\, M
\end{array}.
\label{groupmorpmany}\end{equation}
We call a map $f\colon \Sigma\rightarrow\Sigma$  {\it liftable}\/ to
$\Pi(\Sigma)$ if there exists a map
$ F\colon \Sigma\rightarrow\Pi(\Sigma)$, such that $s\circ
 F=\id$ and $t\circ F=f$. In this case, ${\cal F}([c_{uv}])= F(u)[c_{uv}]
 F(v)^{-1}$ so that there exists a gauge
transformation $\hat{\Phi}={\cal X}\circ F\in \G^\Sigma$ such that
\begin{equation}
{\cal X}\circ {\cal F}([c_{uv}])=
 {\cal X}([ F(u)]) {\cal X}([c_{uv}]) {\cal X}([
 F(v)^{-1}])= {\cal X}_{\hat\Phi}([c_{uv}]).
\label{definame123}\end{equation}

 Let us comment more on the liftability condition of $f\colon \Sigma\rightarrow
\Sigma$ and give some examples.  Consider  $\Sigma=\sphere^2$. Then
 $\Pi(\sphere^2)=\sphere^2\times\sphere^2$ and every map $f \colon  \sphere^2 \rightarrow \sphere^2$ is liftable
 with $F(u)=(f(u),u)$. In particular a constant map $f(u)=u_0$ is liftable
 and thus any groupoid morphism
$(X, {\cal X})$ is gauge equivalent to a trivial morphism
$(X(u_0), {\cal X}([u_0,u_0]))$ through (\ref{groupmorpmany}) and
\eqref{definame123}. We will use this property to derive the
integrability condition \eqref{stateresult}. Indeed, constant maps
are liftable only for $\sphere^2$.
Another class of liftable maps, for any $\Sigma$, is the
subgroup of diffeomorphisms defined by the group of bisections;
namely, for every $F \in \bis (\Pi(\Sigma))$ we have that
$f=t\circ F\in \Diff(\Sigma)$ is liftable by definition. In general,
not every diffeomorphism of $\Sigma$ comes from a bisection of $\Pi(\Sigma)$.
However, every diffeomorphism connected to the identity is indeed liftable
and so it acts as
a finite gauge
transformation. This is consistent with having a topological theory.
Observe that for $\Sigma$ simply connected every diffeomorphism is liftable, so every
diffeomorphism acts as a finite gauge transformation.

 Now we study the value of the action (\ref{definPS}) on a solution
  $\hat{X} \in {\rm Mor}(T\Sigma,T^*M)$.
 The base map $X$ maps $\Sigma$ to a  symplectic
 leaf ${\cal L}$ (see, e.g.,  the Appendix in \cite{Bonechi:2003hd} for a proof).
 The tangent space of a leaf ${\cal L}$ is
 spanned by the Hamiltonian vector fields $X_f = \{f,\,\,\}$ and is endowed with a symplectic form
 defined by
\begin{equation}
 \omega_{\cal L} (X_f , X_g) = \{ f, g\} .
\label{symplformLL}\end{equation}
 Thus the action (\ref{definPS}) evaluated on a solution $\hat X$
  can be rewritten as the pull-back of the  symplectic
 form $\omega_{\cal L}$ on the leaf ${\cal L}\supset X(\Sigma)$ as
\begin{equation}
 \Scl (X,\eta) = \int\limits_{\Sigma} X^* (\omega_{\cal L}).
\label{pullbsympa}\end{equation}
 It is important to note that not every map $X \colon  \Sigma \rightarrow {\cal L} \subset M$
 is the base map of a Lie algebroid morphism, unless $\dim {\cal L}= \dim M$.

 Next we have to require that for two gauge equivalent solutions
 condition (\ref{valclase}) is satisfied.
 Consider two gauge equivalent solutions $(X,\eta),
 (X_{\hat\Phi}, \eta_{\hat\Phi})\in{\rm Mor}(T\Sigma,T^*M)$
 related by a gauge transformation $\hat{\Phi}\in\G^\Sigma$
 corresponding to a liftable map $f\colon \Sigma\rightarrow\Sigma$
  as described in (\ref{groupmorpmany}) and (\ref{definame123}).
 In its turn the action functionals on these solutions are related to each other by
\begin{equation}
\Scl (X_{\hat\Phi},\eta_{\hat\Phi}) = \int\limits_\Sigma
X^*_{\hat\Phi}(\omega_{\cal L}) = \int\limits_\Sigma f^*\circ
X^*(\omega_{\cal L}) = \deg(f) \Scl (X,\eta) ,
\label{largshatsla}\end{equation}
 where $\deg(f)$ is the degree of the map $f$.
 For $|\deg(f)|=0,2$ the requirement (\ref{valclase}) implies the integrality
 condition
\begin{equation}
 \frac{1}{2\pi \hbar} \Scl (X,\eta)\in \Z
\label{neconw9120}\end{equation}
for any Lie algebroid morphism
$(X,\eta)$.
Since for $\Sigma=\sphere^2$ every constant map is liftable, we get at condition \eqref{stateresult}.

This result can be easily extended
  to surfaces with boundary, see the discussion in Section \ref{s:boundary}. Let us anticipate
  the case of the disk $D$ with boundary
conditions (\ref{CFboundary}), which is the relevant case for
(\ref{defpath}). Due to this particular boundary condition,
  every classical solution is  constant on the boundary $\d D$, $X|_{\d D} =x_0$.
 Therefore classical solutions on the disk with the given boundary conditions
also correspond to Lie algebroid morphism  $T\sphere^2 \rightarrow T^*M$.
  Thus  the previous discussion for $\sphere^2$
 can be  applied to this case and we encounter again the same integrality condition
 (\ref{neconw9120}). More precisely the relevant groupoid for these boundary conditions
 is
\begin{equation}\label{gdm}
\G^{D,M}=\{\hat{\Phi}\in\G^D\ |\ \hat{\Phi}(\partial D)\subset
 M\}
\end{equation}
over $M^D$; the group of bisections can be described as follows
\begin{equation}
\label{gauge group} \bis(\G^{D,M}) =\{\sigma\in\bis(\G^D)\ |\
\sigma(X)|_{\partial D}=X|_{\partial D} \}\;.
\end{equation}
It is clear that (\ref{gauge_trans_gpd}) defines an action of
$\G^{D,M}$ on the space of groupoid morphisms $(X,{\cal X})\in{\rm
Mor}(D\times D,\G(M); M)$ that are trivial on the boundary: i.e.,
${\cal X}(\partial D\times\partial D)\subset M$. We refer to
Section \ref{s:boundary} for a detailed discussion of the case
with generic boundary conditions.

In order to conclude that the integrality condition
(\ref{stateresult}) is actually a necessary condition we have to
discuss if $\G^{D,M}$ acts off-shell and is really a symmetry of
the problem. These are the gauge transformations that look most
natural from the geometrical point of view of Lie groupoid theory,
but it is not obvious a priori that we are allowed to use them.
Since the algebra of infinitesimal transformations closes only
on-shell, the first problem one has to face is the correct
understanding of the infinitesimal off-shell transformations, see
for instance the discussion in \cite{Bojowald:2004wu}. Then, as we
already verified on-shell, the integration of the infinitesimal
transformations (\ref{gaugetransf}) does not have a unique answer
in the same way as the integration of a Lie algebra may correspond
to different Lie groups.

Our strategy will be different. We plan to discuss finite
off-shell gauge transformations elsewhere, while in the next
subsections we will discuss the role of the integrality in the
nonperturbative properties of (\ref{defpath}). We will get a
picture that confirms this choice of the gauge transformations for
a disk $D$. However we cannot claim anything about the case of a
closed surface $\Sigma$, in particular $\sphere^2$.

Finally, we refer to \cite{Bonechi:2005} for a discussion of the
moduli space of Lie algebroid morphisms divided by finite gauge
transformations in the case of a generic surface $\Sigma$.

\subsection{The symplectic case}
\label{subs:symplectic}

In  this subsection we consider the case when the Poisson tensor is
invertible. Then $\omega\equiv\alpha^{-1}$ is a symplectic form.
In this case any map $X\colon \sphere^2 \rightarrow M$ defines a
Lie algebroid morphism $\hat{X}=(X,\eta):T\sphere^2\rightarrow
T^*M$ with $\eta_\mu = - \omega_{\mu\nu}dX^\nu$. Therefore the
integrality condition (\ref{stateresult}) corresponds to the
integrality of the pairing between $[\omega/2\pi\hbar]$ and
$\pi_2(M)$ and is weaker than the usual geometric quantization
condition $\left[\frac{1}{2\pi \hbar} \omega\right] \in H^2(M,
{\mathbb Z})$. These two conditions are the same only if $M$ is
simply connected.


Let us interpret this condition in the path integral defined in
(\ref{defpath}).
 Assume for the moment that the target manifold $M$ is simply connected.
 In the symplectic case the formula (\ref{defpath}) essentially reduces to the original
 Feynman path integral formula for quantum mechanics.
 We can formally integrate $\eta$ and arrive to the expression
\begin{equation}
 f * g (x) = \int\limits_{X(\infty)=x} f(X(1)) g(X(0)) e^{\frac{i}{\hbar} \int\limits_D X^*(\omega)} dX ,
\label{defpath12134}\end{equation}
 where one should integrate over any map $X\colon  D \rightarrow M$. While formally
 integrating\footnote{It is convenient to introduce polar coordinates $(r, \phi)$   on the disk $D$ ($r \leq 1$).
 Then $\eta_\phi$ is regarded as a Lagrangian multiplier with the boundary condition $\eta_\phi|_{r=1}=0$.
 Upon the integration of $\eta_\phi$ we get a $\delta$-function in path integral imposing the relation
 on $\d_r X$ and $\eta_r$ on $r < 1$.}
 $\eta$ in (\ref{defpath}) we took into account the boundary conditions (\ref{CFboundary}).
It is clear that the observable we are integrating depends only on
the boundary values $\gamma=X|_{\partial D}$ so that, if $\omega$
is exact then (\ref{defpath12134}) can be defined as
\begin{equation}
 f*g (x) = \int\limits_{\gamma(\infty)=x} f(\gamma(1)) g(\gamma(0)) e^{\frac{i}{\hbar}
 \int\limits_\gamma d^{-1} \omega} d\gamma \;,
\label{symppath}
\end{equation}
where the sum is over all loops $\gamma\colon\sphere^1\rightarrow M$,
with $\gamma(\infty)=x$, or, equivalently, as proposed in
\cite{Cattaneo:1999fm},  over trajectories $ \gamma\colon
\mathbb{R} \rightarrow M$ with $\gamma(\pm \infty)=x$. Integrality
(\ref{stateresult}) allows one to define the weight function also when
$\omega$ is not exact. Indeed we can define the action using the
following prescription
\begin{equation}
 S= \int\limits_{\gamma} (d^{-1}\omega) \equiv \int\limits_{D} X^*(\omega) ,
\label{definact}
\end{equation}
where $X \colon D \rightarrow M$ is any map such that $X|_{\d
D}=\gamma$. The definition (\ref{definact}) is good if it is
independent of the choice of $X$. Namely, let us choose two maps
$X_1$ and $X_2$ such that $X_1|_{\d \Sigma} = X_2|_{\d \Sigma}=
\gamma$ then the difference between two definitions is given by
\begin{equation}
\int\limits_{D} X_1^* (\omega) - \int\limits_{D} X_2^* (\omega)
=
 \int\limits_{X_1(D)} \omega - \int\limits_{X_2(D)} \omega =
  \int\limits_{c_2} \omega ,
 \label{normaliza}\end{equation}
where $c_2=X_1(D)\cup X_2(D)$ is the two-sphere obtained by
joining $X_1(D)$ and $X_2(D)$ along the boundary. In order to
define unambiguously the path integral measure we have to require
\begin{equation}
  e^{\frac{i}{\hbar}  \int\limits_{c_2} \omega} = 1
\label{defpathe}
\end{equation}
for any sphere $c_2$ in $M$.  This is the standard argument which
is used for the Dirac charge quantization of
 $U(1)$ monopole \cite{Dirac:1931kp}.

If $\pi_1(M)$ is nontrivial, not every loop $\gamma$ can be covered
by a disk $D$, so the path integral in (\ref{symppath}) is over
contractible loops only. If $H_1(M)$ is trivial then any loop is
the boundary of some surface $\Sigma$, possibly with handles, so
that one can extend (\ref{symppath}) to any loops by covering them
with higher genus surfaces. In order to do this we have to require
that $\left[\frac{1}{2\pi \hbar} \omega\right] \in H^2(M, {\mathbb
Z})$. However, if $H_1(M)$ is nontrivial, the extension of the path
integral (\ref{symppath}) to all loops goes beyond the Poisson
sigma model. The associativity of Kontsevich formula suggests that
it is not necessary to include these contributions in order to
have an associative $*$-product.

We can make contact with the discussion of the previous subsection
by reinterpreting the equivalence suggested by this construction
in terms of gauge invariance under finite gauge transformations.
Indeed thanks to integrality, the 1d path integral
(\ref{symppath}) can be defined as the 2d path integral
(\ref{defpath12134}), and we can consider this construction as a
way of getting rid of the gauge equivalence. So the correct gauge
transformation is the one that identifies two maps $X_i\colon
D\rightarrow M$, $i=1,2$, if they coincide on the boundary. We
consider again $M$ to be simply connected, so that $\G(M)=M\times
M$ is the (ssc) groupoid integrating it.

It is immediately clear that the homotopy equivalence is too weak.
{}From the previous discussion, one possible candidate for the
gauge group could be $(\bis(M\times M))^D=\{\sigma\colon D
\rightarrow \Diff(M)\,|\,\sigma|_{\d D} = \id\}$, that integrates
the infinitesimal gauge transformations of the model, $\delta
X^\mu = \alpha^{\mu\nu}\beta_\nu$ with $\beta|_{\d D}=0$. In
general also this group will be too small; in fact, take any two
sphere in M, divide it in two disks $X_1(D)$ and $X_2(D)$, we
require that it exists $\sigma\in(\bis(M\times M))^D$ such that
$X_2=\sigma(X_1)$. In particular, let the loop be the point $x$,
we require that any class in $\pi_2(M,x)$ be represented by
$\sigma(x)$. This will be impossible, if for example
$\pi_2(\Diff(M))=0$, like for $M=\sphere^2$, (see \cite{Smale}).
Finally, let us analyze the groupoid $\G^{D,M}$ defined in
(\ref{gauge group}). It clearly acts on $M^D$ and it is easy to
verify that the source fiber $\G^{D,M}_X$ of $X\in M^D$ is
$\{(Y,X) \,|\, Y:D\rightarrow M, Y|_{\d D}= X|_{\d D}\}$ so that
any two configurations $X_1$ and $X_2$ coinciding on the boundary
will be gauge equivalent.

We should stress that this discussion of the finite gauge
transformations is based on the action (\ref{definact}). If we
accept the formal derivation of (\ref{defpath12134}) from
(\ref{defpath}), the instantons of the PSM correspond to those
$X\in M^D$ such that $X|_{\d D}=x$. Therefore they are all
equivalent and are represented by the constant loop $\gamma=x$ in
(\ref{symppath}). This agrees with the discussion of the previous
subsection: the on-shell gauge groupoid $\G^{D,M}$ is extended
off-shell only after the formal integration of $\eta_\phi$ and
there is  only one classical solution with the prescribed boundary
condition, modulo gauge transformations.

\subsection{The general Poisson case}
\label{subs:poisson}

The discussion of (\ref{defpath}) of the previous subsection for
the symplectic case can be successfully repeated for a generic
integrable Poisson manifold.

Let us start from (\ref{defpath}) and introduce polar coordinates
$(r,\phi)$ on the disk $D$; once that we integrate over the
Lagrange multiplier $\eta_\phi$ we get
$$
f*g(x) = \int\limits_{X(\infty)=x} f(X(1))g(X(0)) e^{\frac{i}{\hbar} \int\limits_D d^2\sigma
\partial_\phi X^\mu \eta_{r\mu}} \delta(\partial_r X+\alpha\eta_r)dX~d\eta_r\;.
$$
Although the Poisson tensor is degenerate, the constraint
$\partial_r X+\alpha\eta_r=0$ can be explicitly solved in the
following way. We are going to prove that for any $(X,\eta_r)$
satisfying $X(\infty)=x$ and $\partial_r X+\alpha\eta_r=0$ there
exists a unique algebroid morphism
$(X,\eta_r,\tilde{\eta}_\phi):TD\rightarrow T^*M$ or, equivalently, a map
$\xi_{X,\eta_r}:D\rightarrow \G(M)_x$ satisfying
$\xi_{X,\eta_r}(\infty)=x$, and vice versa.
 It should be stressed that $\tilde{\eta}_\phi$ is determined
 by $\eta_r$ and $X$ and should not be confused with $\eta_\phi$.
Indeed, for each fixed
$\phi\in\sphere^1$, we let $X^\phi(r)=X(re^{i\phi})$,
$\epsilon^\phi(r)=\eta_r(re^{i\phi})$; it is clear that
$(X^\phi,\epsilon^\phi)$ is an algebroid morphism from
$TI\rightarrow T^*M$ and we consider the unique groupoid morphism
$(X^\phi,{\cal X}^\phi):I\times I\rightarrow \G(M)$ that
integrates it. Then we define $\xi\colon\sphere^1\times I\rightarrow
\G(M)_{X(0)}$, $ \xi(\phi,r)={\cal X}^\phi(r,0)$. Since
$\xi(\phi,0)={\cal X}^{\phi}(0,0)=X(0)$, we can consider $\xi\colon
D\rightarrow \G(M)_{X(0)}$, so that it defines a groupoid morphism
$D\times D\rightarrow \G(M): (u,v)\rightarrow \xi(u)\xi(v)^{-1}$
and by differentiation the desired algebroid morphism
$(X,\eta):TD\rightarrow T^*M$, where $X(u)=t(\xi(u))$ and
$\eta|_{T_uD}=R_{\xi(u)^{-1}*}\circ\xi_*$. Finally,
$\xi_{X,\eta_r}:D\rightarrow \G(M)_x $ defined by
$\xi_{X,\eta_r}(u)=\xi(u)\xi(\infty)^{-1}$ satisfies
$\xi_{X,\eta_r}(\infty)=x$. The converse statement is at this
point obvious.  We remark that $\tilde{\eta}_\phi|_{\partial D}$ is not
zero in general, and the on-shell configurations of the PSM
correspond to those algebroid morphisms that satisfy boundary
conditions (\ref{CFboundary}) or, equivalently, to those maps
$\xi$ such that $\xi(\partial D)=x$.

The value of the action on these configurations
$(X,\eta_r)\leftrightarrow\xi$ can be equivalently written as
\begin{equation}
 S(X,\eta_r)= \int\limits_D \partial_\phi X^\mu\eta_{r\mu} rdrd\phi= - \int\limits_D
 \alpha^{\mu\nu} \tilde{\eta}_{\phi\mu}\eta_{r\nu} rdrd\phi = - \int\limits_D
 X^*(\omega_{{\cal L}_x})=\int\limits_D\xi^*(\Omega)=S(\xi),
\end{equation}
where $\omega_{{\cal L}_x}$ is the symplectic form of the leaf
${\cal L}_x$ containing $x$ and $\Omega=t^*(\omega_{{\cal L}_x})$
is the symplectic form of the symplectic groupoid restricted to
$\G(M)_x$.

We can then rewrite (\ref{defpath}), analogously to
(\ref{defpath12134}), as a sum over all algebroid morphisms
$(X,\eta):TD\rightarrow T^*M$, or, equivalently over all $s$-vertical
maps $\xi\colon D\rightarrow \G(M)_x$ with $\xi(\infty)=x$,
\begin{equation}
\label{defpath12134poisson} f*g(x)= \int\limits_{\xi:D\rightarrow
\G(M)_x,~ \xi(\infty)=x} f(t(\xi(1))) g(t(\xi(0)))
e^{\frac{i}{\hbar} S(\xi)} ~d\xi\;.
\end{equation}
It is now easy to see that the integrality condition
(\ref{stateresult}) plays the same role as in the symplectic
case. We recall that (\ref{stateresult}) implies the integrality
the symplectic form $\Omega$ over all the spherical cycles
contained in $\G(M)_x$. If $\alpha$ satisfies (\ref{stateresult})
then $\exp i/\hbar S(\xi)$ depends only on the boundary value
$\xi|_{\partial D}$. In fact if $\xi_1|_{\partial
D}=\xi_2|_{\partial D}$ then
$S(\xi_1)-S(\xi_2)=\int\limits_{\xi_1\cup\xi_2}t^*(\omega_{{\cal L
}_x})\in2\pi \hbar\Z$, where $\xi_1\cup\xi_2$ is the sphere in $\G(M)_x$
obtained by joining $\xi_1$ and $\xi_2$ along the boundary. By
identifying the configurations coinciding on the boundary we get
the final formula
\begin{equation}
\label{defpath12134poisson2} f*g(x)=
\int\limits_{\xi:\sphere^1\rightarrow \G(M)_x,~\xi(\infty)=x}
f(t(\xi(1))) g(t(\xi(0))) e^{\frac{i}{\hbar}S(\xi)} d\xi \;.
\end{equation}
Finally, it is clear that, since classical solutions are constant
on the boundary, they all correspond to the constant map.

It is useful to remark that one can get formula
(\ref{defpath12134poisson2}) by applying the construction of
$\G(M)$ as a Marsden--Weinstein reduction in $T^*PM$, as
explained in Section \ref{s:groupoid}. In fact, after the formal
integration of $\eta_\phi$, we can interpret each configuration
$(X,\eta_r)$ as a loop $\gamma:\sphere^1\rightarrow \C_{X(0)}$,
where $\gamma(\phi)=(X^\phi,\epsilon^\phi)$ and $\C_y$ denotes the
set of algebroid morphisms starting at $y$. Let
$\gamma(\phi)\rightarrow\underline{\gamma(\phi)}$ denote the quotient
$\C\rightarrow\G(M)$. Then the above construction means that we
identify two loops $\gamma_1$ and $\gamma_2$ if and only if
$\underline{\gamma_1(\phi)}\ \underline{\gamma_1(\infty)}^{-1}=
\underline{\gamma_2(\phi)}\ \underline{\gamma_2(\infty)}^{-1}$ for
each $\phi\in \sphere^1$.

We can now analyze this construction in terms of the gauge
transformations introduced in Section \ref{s:gauge}. The groupoid
$\G^{D,M}$ acts on $\xi\colon D\rightarrow \G(M)_x$ as follows:
let $\hat{\Phi}\in\G^{D,M}$ be such that
$s(\hat{\Phi})=t\circ\xi$, then we define
\begin{equation}
\label{action_offshell} \xi_{\hat\Phi}(u) = \hat\Phi(u)\xi(u)
~~~~~u\in D\;.
\end{equation}
We see that this action coincides with (\ref{gauge_trans_gpd})
where ${\cal X}(u,v)=\xi(u)\xi(v)^{-1}$ and $X(u)=t\circ\xi(u)$.
Let us verify that this action induces the equivalence relation
involved in the path integral construction. In fact if
$\xi|_{\partial D}=\nu|_{\partial D}$, then we can write
$\nu=\xi{}_{\hat\Phi}$ where $\hat\Phi\in\G(M)^{D,M}$ is defined
by $\hat\Phi(u)=\nu(u)\xi(u)^{-1}$.

\bigskip\bigskip

\section{Prequantization of symplectic groupoids}
\label{s:groupoid}

In this Section we study the prequantization of the source simply
connected (ssc) symplectic groupoid $\G(M)$ of an integrable
Poisson manifold $M$. Since the space of units is a Lagrangian
submanifold, every prequantization of $\G(M)$ induces a flat
structure on $M$. In \cite{WXu} it has been shown that if $\G(M)$
is prequantizable, then there is a unique prequantization bundle
of $\G(M)$ such that this flat structure is trivial. In
\cite{crainic1} it has been shown that a necessary and sufficient
condition for the symplectic groupoid ${\cal G}(M)$ to be
prequantizable is that the symplectic form is integer on every
$\sphere^2\subset s^{-1}(y)$, for every $y\in M$, where $s$ is the
source map. Since any two-sphere in the source fiber defines an
algebroid morphism and vice versa, this integrality condition is
the same as \eqref{stateresult}.

We develop here an alternative and straightforward approach based
on the description of $\G(M)$ as the Marsden--Weinstein symplectic
reduction of a cotangent bundle, obtained in \cite{sympgroup}. We
are going to get the reduction of the prequantization of the
cotangent bundle. This will give us an explicit description of the
prequantization of $\G(M)$ and, hopefully, will clarify some
important issues.

Let us first recall the construction of the symplectic groupoid
${\cal G}(M)$ in \cite{sympgroup}. Let us consider the Hamiltonian
formulation of the Poisson sigma model over a rectangle
$R=[-T,T]\times I$ with coordinates $(t,u)$ labeling time and
space respectively. The boundary conditions are $\eta_t=0$ on
$[-T,T]\times \d I$. In this world-sheet geometry the action
(\ref{definPS}) can be rewritten as
 \begin{equation}
  S= \int\limits_{R} dt\, du [-\langle \eta_{u}, \d_t  X \rangle
 + \langle \eta_{t}, (\d_u X + \alpha \eta_{u}) \rangle ] .
 \label{hamPosact}\end{equation}
The corresponding phase space is given by the space of vector
bundle morphisms $TI \rightarrow T^*M$ with the symplectic
structure defined by the by first term in the action
(\ref{hamPosact}). This phase space can be interpreted as the
cotangent bundle $T^*PM$ of the path space $PM$ of maps
$I\rightarrow M$ with exact symplectic structure $d\theta$, where
$\theta$ is the canonical Liouville form. The field $\eta_t$ is a
Lagrangian multiplier which enforces the first class constraints
\begin{equation}
 \d_u X + \alpha (X)\eta_{u}  = 0
\label{firstcxka}\end{equation} that generate a distribution of
subspaces spanned by Hamiltonian vector fields. Let $H$ be the
group generated by the Hamiltonian vector fields. The constraints
(\ref{firstcxka}) define a coisotropic infinite dimensional
submanifold ${\cal C}(M)$ of $T^*PM$. The reduced phase space
$\underline{{\cal C}(M)}$ of the Poisson sigma model is the set of
integral manifolds of this distribution and is obtained as the
Marsden--Weinstein symplectic quotient $T^*PM//H=\C(M)/H$. In
\cite{sympgroup} it has been shown that $\underline{{\cal C}(M)}$
carries the structure of topological groupoid over $M$ and, if smooth, it is
the symplectic groupoid ${\cal G}(M)$. Alternatively, using the
language of Lie algebroids one can give the following description
of $\G(M)$ \cite{severa, cranfern}. Elements of ${\cal C}(M)$ are
those bundle maps that are also Lie algebroid morphisms
$\hat{\gamma}\colon  TI\rightarrow T^*M$. Elements of $\underline{{\cal
C}(M)}$ are Lie algebroid morphisms modulo homotopy. We say that
two algebroid morphisms $\hat{\gamma}_0, \hat{\gamma}_1 \colon  TI
\rightarrow T^*M$ are homotopic if there exists a Lie algebroid
morphism $\hat{\gamma}\colon T([0,1]\times I)\rightarrow T^*M$ such that
$\hat{\gamma}|_{\{0\}\times I}=\hat{\gamma}_0$,
$\hat{\gamma}|_{\{1\}\times I}=\hat{\gamma}_1$ and
$\hat{\gamma}|_{[0,1]\times \d I}$ is the zero bundle map (i.e.,
the boundary conditions (\ref{CFboundary}) are satisfied).

Let us recall the relation between symplectic reduction and
prequantization in the finite dimensional case following
\cite{Woodhouse}. Consider the cotangent bundle $T^*{\cal M}$ of a
manifold $\cal M$ with the canonical symplectic structure $\omega
= d \theta$; let $\C$ be a coisotropic submanifold of $T^*\M$ and
$\underline{\cal C}$ be the reduced phase space with symplectic
form $\underline{\omega}$. The cotangent bundle  $T^* {\cal M}$ is
prequantized by the trivial line bundle $T^*{\cal M}\times
{\mathbb C}$ with connection $\nabla^\theta\equiv d-(i/\hbar)\theta$. Then a
sufficient condition for $\underline{\cal C}$ to be prequantizable
is that $\theta$ satisfies the integrality condition
\begin{equation}
 \frac{1}{2\pi \hbar} \int\limits_l \theta \in \mathbb{Z}
\label{inegraleaf}
\end{equation}
for every closed curve $l$ inside a leaf of the characteristic
foliation. In fact we can construct the prequantization line
bundle $L({\underline {\cal C}})$ as ${\cal C}\times {\mathbb
C}/\sim$, where $(c_0, z_0)\sim(c_1, z_1)$ whenever $c_0$ and
$c_1$ are in the same leaf and \begin{equation}
 z_1 = z_0\,\, e^{\frac{i}{\hbar} \int\limits_{c_0}^{c_1} \theta},
\label{equivalelsao}
\end{equation}
where $\theta$ is integrated over any path contained in the leaf
and connecting $c_0$ to $c_1$. The integrality (\ref{inegraleaf})
ensures that the equivalence relation is well-defined, i.e., that the
exponent in (\ref{equivalelsao}) does not depend on the concrete
choice of a path between $c_0$ and $c_1$. A section
$\psi\in\Gamma(L(\underline{\cal C}))$ can be identified with an
equivariant function $\psi\colon {\cal C}\rightarrow{\mathbb C}$,
i.e., a solution to $\nabla^\theta_V\psi=0$ for every vector field $V$ along
the leaves. If $W$ is a vector field on $\cal C$ that projects to
a vector field $W'$ on $\underline{\cal C}$, then $\nabla^\theta_W$
depends only on $W'$ and defines a connection on
$\Gamma(L(\underline{\cal C}))$ whose curvature is
$\underline{\omega}$. For further details we refer to Proposition
(8.4.9) in \cite{Woodhouse}. Although this Proposition is proved
for the finite dimensional case, its generalization for the
infinite dimensional cotangent bundle is straightforward.

In the following Proposition, we apply the above construction to
our case of $\M=PM$, $\C=\C(M)$ and $\underline{\C}=\G(M)$.

\begin{proposition} \begin{itemize}
\item[$i$)] The ssc symplectic groupoid $\G(M)$ is prequantizable if
and only if $\alpha$ satisfies the integrality condition
(\ref{newcond123}) for $c_2=\sphere^2$.
\item[$ii$)] The reduction of the prequantization of $T^*PM$ induces a trivial flat structure over $M$.\end{itemize}
\end{proposition}
We give an alternative proof to the one in \cite{crainic1}.

{\it Proof.}  Let $\alpha$ satisfy (\ref{newcond123}). The
coisotropic leaves of ${\cal C}(M)$ are homotopic Lie algebroid
morphisms $\hat{\gamma}\colon T I \rightarrow T^*M$, so that a loop
contained in a leaf is a Lie algebroid morphism $T (S^1\times I)
\rightarrow T^*M$. Due to the boundary conditions, the boundaries
of the annulus $S^1\times I$ are mapped to fixed points, so we
effectively deal with a Lie algebroid morphism $\hat{X} = (X,
\eta)\colon  T\sphere^2 \rightarrow T^*M$. The integral of the Liouville
form $\theta$ is given by the first term of (\ref{hamPosact}) and
thus on a leaf it coincides with the pairing between $\alpha$ and
$\hat{X}$. Therefore the condition (\ref{inegraleaf}) becomes
(\ref{newcond123}), i.e.,
  \begin{equation}
\frac{1}{2\pi \hbar} \int\limits_l \theta = \frac{1}{2\pi\hbar}
\Scl  (X, \eta)= \frac{1}{2\pi\hbar} \ll \alpha, (\sphere^2,
\hat{X})   \gg   \,\,\,
 \in\,\,\, \mathbb{Z}   ,
  \label{mainress4}\end{equation}
so that $\G(M)$ is prequantizable.

Conversely, let us suppose that $\G(M)$ is prequantizable with
integer symplectic form $\Omega$. Let $\C_y\subset\C(M)$ be the
subset of those algebroid morphisms starting at $y\in M$ and let
$\gamma\colon \sphere^1\rightarrow \C_y$ be defined by
$\gamma(t)=(X(t,u),\eta(t,u))$. Then
$F(r,t)=(X(t,ru),r\eta(t,ru))$ is a homotopy to the trivial
morphism $(y,0)$, {i.e.}, $\C_y$ is simply connected for every
$y\in M$. Since homotopic morphisms have the same end points,
every leaf is contained in $\C_y$ for some $y\in M$ and every loop
$l$ contained in it can be covered by a disk $D\subset \C_y$, such
that $\partial D=l$. We then have that $ \int\limits_l \theta =
\int\limits_D d\theta = \int\limits_{\sphere^2} \Omega \in 2\pi \hbar \Z$, so that $\theta$
satisfies (\ref{inegraleaf}).

Finally, it can be verified that $L(\G(M))|_M=M\times {\mathbb C
}$ and $\nabla^\theta\psi=d\psi$, for every $\psi\colon M\rightarrow{\mathbb C
}$; {i.e.}, the flat connection is trivial. $\square$

\begin{remark}{\rm
We can explicitly describe the line bundle $L(\G(M))$ as ${\cal
C}(M)\times{\mathbb C}/ \sim$, where $(\hat{\gamma}_0, z_0) \sim
(\hat{\gamma}_1, z_1)$ if $\hat{\gamma}_0$ and $\hat{\gamma}_1$
are homotopic algebroid morphisms and
     \begin{equation}
       z_1 = z_0 e^{\frac{i}{\hbar}   \ll \alpha, (D, \hat{X}),
   \gg   }
     \label{eqdjcotnaldf}\end{equation}
where $\hat{X}\colon TD \rightarrow T^*M$ is the algebroid homotopy
between $\hat{\gamma_0}$ and $\hat{\gamma_1}$ seen as an algebroid
morphism from the disk thanks to the boundary conditions.  Because of
(\ref{mainress4}) this definition of equivalence does not depend
on the concrete choice of $(D, \hat{X})$.}
\end{remark}

Furthermore we can reinterpret the integrality condition
(\ref{inegraleaf}) using the following path integral argument.
Consider a cotangent bundle $T^*{\cal M}$ with the canonical
symplectic form $\omega=d\theta$. We assume that the coisotropic
submanifold $\C$ is defined by the first class constraints
$\{\phi_i\}$ and let $G$ be the group of diffeomorphisms generated
by the constraints. Let $p\colon \C\rightarrow\underline{\C}=\C/G$ be
the canonical projection and let $\underline\omega$ be the
symplectic form of $\underline\C$ ({i.e.},
$p^*\underline\omega=\omega$).

We consider the following action functional \begin{equation}
 S(\gamma, \lambda)= \int\limits_{\gamma}( p \dot{q} - \lambda^i \phi_i (q, p))dt \;.
\label{theory12}\end{equation}
The ultimate goal is to define the path integral
 \begin{equation}
  \int  d\gamma\,\, d\lambda\,\, \cdots \,\,e^{\frac{i}{\hbar}S(\gamma,\lambda)} =
 \int d\gamma\,\, \delta(\phi^i(q,p)=0)\,\,\cdots\,\, e^{\frac{i}{\hbar} S(\gamma, 0)},
 \label{defpathintex}\end{equation}
where $\cdots$ stands for some gauge invariant observable. The
integration over the $\lambda_i$'s restricts the path integral to
${\cal E} = \{ \gamma\colon  I \rightarrow {\cal C}\}$; the gauge
symmetries force us to identify $\gamma(t)$ and $g(t)\gamma(t)$,
with $g\in PG=\{g\colon I\rightarrow G\}$ and so to define the integral
over ${\cal E}/PG$, which are the maps with values in the reduced
phase space $\underline{\cal C}$. A necessary condition for this to
work is that the procedure defines a correct measure over ${\cal
E}/PG$.

Let us study in details the measure of the path integral
(\ref{defpathintex}). Once the integration over $\lambda$ is done,
the action can be defined as follows:
\begin{equation} S(\gamma) \equiv  S
(\gamma, 0)= \int\limits_{\gamma} \theta ,
\label{difdefacex}\end{equation}

We need to compare the actions of two gauge equivalent
configurations $\gamma$ and $g\gamma$, $g \in PG$. Namely we have
\begin{equation}
 \int\limits_\gamma \theta - \int\limits_{p_0} \theta - \int\limits_{g\gamma}\theta +
 \int\limits_{p_1} \theta =
   \int\limits_\Sigma \omega = \int\limits_{i(\Sigma)}  \underline{\omega} = 0
\label{defiaod[20}
\end{equation}
where $p_0$ ($p_1$) is a path contained in a leaf and connecting
$\gamma(0)$ ($\gamma(1)$) and $g(0)\gamma(0)$ ($g(1)\gamma(1)$);
$\Sigma$ is a surface such that $\d\Sigma = \gamma
p_0^{-1}(g\gamma)^{-1}p_1$. Such $\Sigma$ exists since $G$ is
connected so that $g$ is homotopic to identity. From
(\ref{defiaod[20}) we have \begin{equation}
 e^{\frac{i}{\hbar} S(g\gamma)} = e^{\frac{i}{\hbar}\int\limits_{p_1} \theta}
 e^{-\frac{i}{\hbar}\int\limits_{p_0} \theta} e^{\frac{i}{\hbar}S(\gamma)},
\label{cond2937}
\end{equation}
and if the condition (\ref{inegraleaf}) is satisfied then the
expression (\ref{cond2937}) does not depend on the concrete choice
of paths $p_0$ and $p_1$, only end points matter. Remark that
fixing the end points of $\gamma$ in the path integral would break
gauge invariance. We can now do the integration along the fiber and
reduce the integral to ${\cal E}/PG$; {i.e.}, we can define the
measure
\begin{equation}
 W(\underline{\gamma}) = e^{\frac{i}{\hbar} S(\gamma)}
\int\limits_{G\times G} dg_0dg_1
e^{-\frac{i}{\hbar}\int\limits_{\gamma(0)}^{\gamma(0)g_0}\theta}
e^{\frac{i}{\hbar}\int\limits_{\gamma(1)}^{\gamma(1)g_1}\theta}
  \;.
\label{mesianero}
\end{equation}

Vice versa in order to separate integrations along and transverse to
a fiber in (\ref{inegraleaf}) one needs the integrality condition.

The present argument based on path integral is not a rigorous one.
However, in our view, it offers a nice physical intuition behind the
integrality condition (\ref{inegraleaf}). Although we stated this
for a quantum mechanical system (i.e., with finite dimensional
$T^*{\cal M}$), it is straightforward to generalize it for a
Poisson sigma model where the cotangent bundle and the coisotropic
submanifold are infinite dimensional.

\section{D-branes and integrality}
\label{s:boundary}

In the previous Sections we discussed the Poisson sigma model over $\Sigma$
 either without boundary or with the specific boundary conditions given in (\ref{CFboundary}).
 However, there exist more general non-symmetry-breaking boundary conditions, as discussed in
\cite{Cattaneo:2003dp}, which correspond to coisotropic submanifolds of $M$. A submanifold ${\cal D}$
 of $M$ is called coisotropic if $\sharp N^*{\cal D} \subset T{\cal D}$, where $N^*{\cal D}$ is the conormal
 bundle of ${\cal D}$. The boundary conditions corresponding to the coisotropic submanifold
 ${\cal D}$ are
\begin{equation} X(\d\Sigma) \subset {\cal D},\,\,\,\,\,\,\,\,\,\,\,
\eta_t|_{\d \Sigma} \in \Gamma(X^*(N^*{\cal
D})),\,\,\,\,\,\,\,\,\,\,\, \beta|_{\d\Sigma} \in
\Gamma(X^*(N^*{\cal D})),
\label{coistbound}\end{equation}
 and (\ref{CFboundary}) corresponds to the case when ${\cal D}= M$. Observe that
 $\sharp N^*{\cal D}$ defines an involutive distribution on ${\cal D}$ (the characteristic
 foliation) and the gauge transformations correspond to diffeomorphisms along the leaves.

 The results of
both Sections \ref{s:gauge} and \ref{s:groupoid} can be generalized for
 the case of a generic coisotropic brane.
 In following two subsections we sketch the construction.

\subsection{On-shell gauge transformations with branes}
\label{subs:gaugebrane}

 In analogy with the discussion in subsection \ref{subs:onshell} we define
 an on-shell gauge group and deduce an integrality condition for
 the classical action.

 Consider a surface $\Sigma$ with $n$ boundary components,
$\partial\Sigma=\bigcup\limits_{i=1}^n \sphere^1_i$ and associate to $\d \Sigma$
a set   $\C=\{C_i\}$ of $n$ coisotropic submanifolds of $M$.
 The conormal subbundles $N^*C_i$ 
 are Lie subalgebroids of $T^*M$ with
anchor map $\sharp\colon N^*C_i\rightarrow TC_i$. Assuming the Poisson manifold
 $M$ to be integrable, there are Lagrangian Lie subgroupoids
$\G(C_i)\subset\G(M)$ that integrate the Lie subalgebroids $N^*C_i$ \cite{Cattaneo:2004}.
 As we said, the boundary conditions are defined by a choice
of $\C$, {i.e.}, $X\colon \sphere^1_i\rightarrow C_i$ and
$\eta|_{\sphere_i}\in\Gamma(T^*\sphere^1_i\otimes X^*(N^*C_i))$. A classical
solution of the equations of motion with these boundary conditions
 is given by a Lie algebroid morphism $\hat{X}=(X,\eta)\colon T\Sigma\rightarrow T^*M$
such that its restriction on the $i$-th component of the boundary
$(X,\eta)|_{\sphere^1_i}$ is a Lie algebroid morphism from $T\sphere^1_i$ to
$N^*C_i$.  We denote the space of such Lie algebroid morphisms by ${\rm
Mor}(T\Sigma,T^*M; N^*\C)$.
 As a consequence of the integrability of the Poisson manifold $M$, we have
that this space is the same as the space ${\rm
Mor}(\Pi(\Sigma),\G(M); \G(\C))$ of groupoid morphisms $(X,{\cal X})$
such that
${\cal X}\left(\Pi(\Sigma)|_{\sphere^1_i}\right)\subset\G(C_i)$, where
$\Pi(\Sigma)|_{\sphere^1_i}=\{[c_{uv}]\in\Pi(\Sigma)\,|\,
c_{uv}\subset\sphere^1_i\}$ is the subgroupoid of $\Pi(\Sigma)$
integrating the $i$-th boundary component $\sphere^1_i$.

In analogy with the closed case, we define the groupoid
$\G^{\Sigma,\C}=\{\hat{\Phi}\colon \Sigma\rightarrow\G(M)\,|\,
\hat{\Phi}(\sphere^1_i)\subset\G(C_i)\ \forall i\}$ over
$M^{\Sigma,\C}=\{\Phi\colon \Sigma\rightarrow M\,|\,
\Phi(\sphere^1_i)\subset C_i\ \forall i\}$. Take $(X,{\cal
X})\in{\rm Mor}(\Pi(\Sigma),\G(M);\G(\C))$ and
${\hat\Phi}\in\G^{\Sigma,\C}$ with $s(\hat\Phi)=X$, then the
groupoid action of $\G^{\Sigma,\C}$ on ${\rm
Mor}(\Pi(\Sigma),\G(M);\G(\C))$ is defined by
\begin{equation} {\cal X}_{\hat\Phi}([c_{uv}]) =
\hat\Phi(u){\cal X}([c_{uv}])\hat\Phi(v)^{-1} \; .
\label{neq29101}\end{equation}

As in the closed case, we can describe part of the gauge groupoid
through liftable maps.
 A map $f\colon \Sigma\rightarrow\Sigma$ preserving the boundary ({i.e.},
  $f(\sphere^1_i)\subset\sphere^1_i$) is said to be
liftable to $\Pi(\Sigma)$ if there exists a map $F\colon
\Sigma\rightarrow\Pi(\Sigma)$ such that
$F(\sphere^1_i)\subset\Pi(\Sigma)|_{\sphere^1_i}$, $s\circ F=\id$
and $t\circ F=f$. Then ${\cal F}([c_{uv}])= F(u)[c_{uv}]
F(v)^{-1}$ is the corresponding groupoid morphism and it is easy
to verify that $(X\circ f, {\cal X}\circ {\cal F})\in{\rm
Mor}(\Pi(\Sigma),\G(M);\G(\C))$. Moreover, there exists
$\hat\Phi\in\G^{\Sigma,\C}$ such that ${\cal X}\circ{\cal F}={\cal
X}_{\hat\Phi}$ with $\hat\Phi= {\cal X}\circ{\cal F}$.

Let us concentrate on the case of the disk $D$ with boundary
condition defined by the coisotropic submanifold $C$.  The corresponding
 ssc groupoid is $\Pi(D)=D\times D$ and any map $f\colon D\rightarrow D$ preserving
the boundary is liftable to $D\times D$.  Considering a
constant map $f(u)=u_0\in\partial D$ we conclude that any
groupoid morphism $(X,{\cal X})$ is gauge equivalent to the trivial
one $(X(u)=X(u_0), {\cal X}([c_{uv}])={\cal X}([u_0, u_0]))$. As before we
  can use the argument (\ref{largshatsla})
to deduce the necessary integrality condition
\begin{equation} \frac{1}{2\pi\hbar} \Scl (X,\eta)\,\,\,\in\,\,\,\Z
\label{inegralconbound}\end{equation}
for every $(X,\eta)\in{\rm
Mor}(TD,T^*M; N^*\C )$.
 Indeed this integrality condition has cohomological interpretation since it depends
  only on a class $[(\alpha, 0)] \in H_{PL}^\bullet(M, {\cal D}, \alpha)$ in the appropriate relative
   cohomology, see the Appendix for the construction. Using the pairing between
    relative (co)homology groups introduced in
    the Appendix we can write the integrality condition as
    \begin{equation}
   \frac{1}{2\pi\hbar} \ll (\alpha, 0), (D, \hat{X}, \d D, \hat{X}|_{\d D}) \gg\,\,\, \in\,\,\, \Z ,
    \label{inegal9322}\end{equation}
    where $\hat{X}\colon  TD \rightarrow T^*M$ and $\hat{X}|_{\d D}\colon  \d D \rightarrow N^* {\cal D}$
     are Lie algebroid morphisms.
     If $(\alpha, 0)$ is the trivial element in $H_{PL}^\bullet(M, {\cal D}, \alpha)$, then
      the classical action is zero.

Condition (\ref{inegal9322}) is a generalization of condition (\ref{mainress4}) (or condition (\ref{newcond123})
 for $\Sigma =\sphere^2$). In the case ${\cal D} =
 M$ the condition (\ref{inegal9322}) is equivalent to (\ref{mainress4}).

Finally, if $M$ is symplectic with exact symplectic form $\omega=d\theta$, then condition \eqref{inegal9322}
is equivalent to  Woodhouse's condition \eqref{inegraleaf}.

\subsection{Reduced phase space picture}
\label{subs:redbrane}

One can consider the reduced phase picture for the case of generic coisotropic branes.
 We can apply the logic from Section \ref{s:groupoid} to the case of general boundary conditions.
 Namely, consider the Poisson sigma model over rectangle  $R=[-T, T] \times I$
  and impose on the boundary $\d I =\{0,1\}$
  the boundary conditions (\ref{coistbound}) corresponding to the coisotropic
  submanifolds ${\cal D}_0$ and ${\cal D}_1$ respectively.
 Following \cite{Cattaneo:2003dp} we consider the reduced phase space with these boundary conditions.
We denote by    ${\cal C}(M; {\cal D}_0, {\cal D}_1)$ the submanifold of ${\cal C}(M)$ where the base maps are paths connecting
     ${\cal D}_0$ to ${\cal D}_1$. It turns out that ${\cal C}(M; {\cal D}_0, {\cal D}_1)$ is coisotropic in $T^*PM$.
    The characteristic foliation on ${\cal C}(M; {\cal D}_0, {\cal D}_1)$ may
    move the end points of the base maps only along the characteristic
  foliations of ${\cal D}_0$ and ${\cal D}_1$. The leaf space $\underline{{\cal C}(M; {\cal D}_0, {\cal D}_1)}$
  corresponds to a reduced phase space of the model on the annulus with the prescribed boundary conditions.
   If smooth, $\underline{{\cal C}(M; {\cal D}_0, {\cal D}_1)}$ is endowed with a s symplectic structure.
 Using the logic from  Section \ref{s:groupoid}  we can study the condition (\ref{inegraleaf})
  in the present setup and construct the prequantization line bundle of
  $\underline{{\cal C}(M; {\cal D}_0, {\cal D}_1)}$.

In particular consider the case ${\cal D}_1=M$. Then the space
$\underline{{\cal C}(M; {\cal D}_0, M)}$ (if smooth) is prequantizable
 if the integrality condition (\ref{inegal9322}) is satisfied. The proof is completely analogous
  to the one presented in Section \ref{s:groupoid}.


\section{Concluding remarks}
\label{s:end}

The integrality condition (\ref{stateresult}) appears to be a
generalization of different integrality conditions within the
geometric quantization program. In particular it points toward the
program of quantizing Poisson manifolds through their
corresponding symplectic groupoids. In fact, condition
(\ref{stateresult}) is a necessary and sufficient condition for
the symplectic groupoid to be prequantizable. On the other hand,
using the Poisson sigma model as a tool, in this paper we showed
that (\ref{stateresult}) plays different roles in the program of
giving a nonperturbative meaning to the Kontsevich formula.

These two facts are strictly related and we believe that they will
contribute in giving a unified description of the quantization
program.  In fact, within quantum field theory, it is a textbook
fact that the path integral of a gauge invariant theory should be
properly defined over its reduced phase space. While for a generic
field theory it is very hard to find an explicit description of
it, for the Poisson sigma model there is the nice geometric
interpretation of reduced phase space as the symplectic groupoid.
We believe that the role of the integrality condition
(\ref{stateresult}) is to reduce the problem of the
nonperturbative definition of (\ref{defpath}) to that of defining
a quantum mechanical path integral for the symplectic groupoid and
at the ultimate end in quantizing it. 

Many interesting problems arise. First of all, it will be
interesting to analyze the role of integrality in the BV action
introduced in \cite{Cattaneo:1999fm}, and hopefully to confirm the
picture provided in this paper on these more solid grounds.

Second, it will be interesting to understand if the reduced phase
picture can also help in understanding more complicated
observables and in particular on closed surfaces.

\bigskip

\bigskip

{\bf Acknowledgments}:
 We are grateful to Domenico Seminara, Marco Zambon, Chenchang Zhu, Thomas Strobl and Alexei Kotov for
 useful discussions.
 MZ thanks INFN Sezione di Firenze and Zurich University
 where part of this work was carried out. FB thanks Zurich
 University.
The research of MZ was supported by EU-grant MEIF-CT-2004-500267.
AC acknowledges partial support of SNF Grant No.~200020-107444/1 and the IHES for hospitality.

\appendix
\Section{Appendix}

Let $c_k=[a_0,\ldots,a_k]$ be a $k$-simplex in ${\mathbb R}^n$ and
let $G$ be an abelian group (for us $G$ will be either $\R$ or $\Z$). Consider the space
\begin{equation}
C_k(M;G)=\left\{\sum\limits_{g_\alpha\in G}g_\alpha
(c_{k\alpha},\hat{X}_\alpha)\right\} \label{fromalsum}\end{equation}
 of finite formal
combinations of Lie algebroid morphisms
$\hat{X}_\alpha\colon  Tc_{k\alpha}\rightarrow T^*M$ with values in $G$ with
boundary operator
\begin{equation}
\partial(c_k, \hat{X})=(\partial c_k, \hat{X}|_{\partial c_k}).
\label{defboudop}\end{equation}
Observe that $\hat{X}|_{\d c_k}$ is automatically a Lie algebroid morphism
  $T(\d c_k) \rightarrow T^*M$ and that the extension of $\d$
 to all $C_k$ by linearity obviously satisfies $\partial^2=0$.
 Thus we associate a homology $H_\bullet(M;G)$ to the Poisson manifold $M$.
 The definition of the cohomology $H^\bullet(M;G)$ with the adjoint coboundary operator
 is straightforward. Let us show that $H^\bullet(M,{\mathbb R})=H_{LP}^\bullet(M, \alpha)$,
 the usual Lichnerowicz--Poisson cohomology.

Any Lie algebroid morphism $
 \hat{X}=(X,\eta)\colon  T c_k\,\,\rightarrow\,\,T^*M$ induces the map
\begin{equation}\label{onsections}
\hat{X}^*\colon  \Gamma(\wedge^{\bullet}
TM)\,\,\rightarrow\,\,\Gamma(\wedge^{\bullet} T^* c_k) \;
\end{equation}
 and thus we can define the following pairing between $A \in
\Gamma(\wedge^k TM)$ and the chain $(c_k,\hat{X})$
\begin{equation}
\label{integralsigma}
 \ll A, (c_k,\hat{X})\gg \equiv \int\limits_{c_k} \hat{X}^*(A)
 = \int\limits_{c_k} A^{\mu_1\ldots \mu_k}\eta_{\mu_1}\wedge \ldots\wedge \eta_{\mu_k}
 \;,
\end{equation}
 where $A$ can be seen as a real valued $k$-cochain; the
following Lemma proves that the Poisson differential
$[\alpha,\cdot]_s$ is the adjoint coboundary to $\partial$.

\medskip
\begin{lemma}
Let $A\in\Gamma(\wedge^{k-1}TM)$ and $\hat{X}=(X,\eta)\colon  T
 c_k\rightarrow T^*M$ be a Lie algebroid morphism. Then we have
\begin{equation}
\label{poisson_stokes} \ll [\alpha,A],(c_k,\hat{X})\gg =
\ll A,\partial(c_k,\hat{X})\gg \;.
\end{equation}
\end{lemma}

{\it Proof}. Starting from the coordinate expression of the
Schouten bracket it is possible to compute that
\begin{equation}
\label{poisson_stokes_2}
[\alpha,A]^{\mu_1\ldots
\mu_k}\eta_{\mu_1}\wedge\ldots\wedge \eta_{\mu_k}  = - d(A^{\mu_1\ldots
\mu_{k-1}}\eta_{\mu_1}\wedge\ldots \wedge \eta_{\mu_{k-1}})\;,
\end{equation}
 where one should use the equations (\ref{eqmotion}), i.e. the fact that $\hat{X}=(X,\eta)$ is
 a Lie algebroid morphism $Tc_k \rightarrow T^*M$.

Indeed this construction works for any Lie algebroid $E$ over $M$ and (\ref{poisson_stokes}) is an analog
 of Stockes' theorem for a Lie algebroid, where $[\alpha,\cdot]_s$ is replaced by the Lie algebroid
 differential $\delta_E$ and $A\in\Gamma(\wedge^{k-1} E^*)$. In fact every morphism of Lie algebroids
 $\Phi\colon E_1\to E_2$ induces a chain map
 $\Phi^*\colon(\Gamma(\Lambda^\bullet E_2^*),\delta_{E_2})\to(\Gamma(\Lambda^\bullet E_1^*),\delta_{E_1})$ and the
 corresponding map $[\Phi^*]$ in cohomology.
 So in particular a morphism $\Phi\colon Tc\to E$ induces a map
 $[\Phi^*]\colon H^\bullet_{\delta_E}\to H^\bullet_{\mathrm{de\ Rham}}(c)$.

 If we fix a submanifold $i\colon  {\cal D} \hookrightarrow M$, then there is a notion of relative
  (co)homology. A $k$-chain in $({\cal D}, M)$  is a linear combination of pairs
  $(c_k, \sigma_{k-1})$  where
   $c_k$ is a $k$-chain in $M$ and $\sigma_{k-1}$ is a $(k-1)$-chain
   in ${\cal D}$. The boundary operator is defined by
    $\d (c_k, \sigma_{k-1}) = (\d c_k + (-1)^k \sigma_{k-1}, \d \sigma_{k-1})$.
   Consider also pairs $(q_k, w_{k-1})$ of forms, $q_k \in \Omega^k(M)$ and
   $w_{k-1} \in \Omega^{k-1}({\cal D})$ with coboundary operator $d(q_k, w_{k-1})=
    (dq_k, dw_{k-1} - (-1)^ki^* q_k)$. There exists a natural pairing
     \begin{equation}
      \langle (q_k, w_{k-1}), (c_k, \sigma_{k-1})  \rangle
      = \int\limits_{c_k} q_k + \int\limits_{\sigma_{k-1}}
       w_{k-1}
     \label{naturalapeuska}\end{equation}
   with the property
   \begin{equation}
     \langle d(q_k, w_{k-1}), (c_{k+1}, \sigma_{k})  \rangle =
      \langle (q_k, w_{k-1}), \d (c_k, \sigma_{k-1})  \rangle .
   \label{prosstusdi}\end{equation}
   This construction gives rise to relative (co)homology groups $H_\bullet(M, {\cal D}, \mathbb{R})$
    and $H^\bullet(M,{\cal D}, \mathbb{R})$. There are also integer versions of these groups.

  Indeed the above construction can also be generalized to the case of a Lie algebroid with a fixed
   subalgebroid.  Let us sketch the construction.
   We are interested in the Lie algebroid $T^*M$ of a Poisson manifold $(M, \alpha)$ and to
    subalgebroids corresponding to coisotropic submanifolds. For each coisotropic
    submanifold $\cal D$ of $M$, we consider
     the Lie algebroid morphism
\begin{equation}
\begin{array}{lll}
 N^*{\cal D} & \stackrel{\hat{\imath}}{\longrightarrow} & T^*M \\
\,\,\,\,\,\Big\downarrow && \Big\downarrow \\
\,\,\,\,\, {\cal D} & \stackrel{i}{\longrightarrow} & M
\end{array},
\label{Alielagmorsj}\end{equation}
 where $\hat{\imath}$ and $i$ are injective immersions. We consider 
 linear combinations of quadruples
  $(c_k, \hat{X}, \sigma_{k-1}, \hat{Y})$ where $\hat{X}\colon  Tc_k \rightarrow T^*M$
   and $\hat Y\colon T\sigma_{k-1} \rightarrow N^*{\cal D}$ are Lie algebroid morphisms.
    We define  a boundary operator by
 \begin{equation}
 \d (c_k, \hat{X}, \sigma_{k-1}, \hat{Y}) = (\d c_{k}, \hat{X}|_{\d c_k} , \d\sigma_{k-1}, \hat{Y}|_{\d \sigma_{k-1}})
 +(-1)^k(\sigma_{k-1},\hat\imath\circ\hat Y,0,0)
   \label{boundrelposp}\end{equation}
    and arrive to a relative homology $H_\bullet(M, {\cal D})$ associated to $(N^*{\cal D}, T^*M)$.
     A $\d$-closed chain $(c_k, \hat{X}, \sigma_{k-1}, \hat{Y})$ corresponds to the case
      when $\d c_k = (-1)^{k-1}\sigma_{k-1}$
      and $\hat{X}|_{\d c_k} = \hat{Y}$.
      We also consider linear combinations of
   pairs $(A_k, v_{k-1})$, $A_k \in \Gamma (\wedge^k TM)$
    and $v_{k-1} \in \Gamma(\wedge^{k-1} N{\cal D})$ and define a coboundary operator by
    \begin{equation}
    \delta (A_k, v_{k-1}) = ([\alpha, A_k], \delta_{N^*{\cal D}}v_{k-1} - (-1)^k \hat\imath^*A_k).
    \label{diffadla;}\end{equation}
This defines the relative Lichnerowicz--Poisson cohomology
$H^\bullet_{LP} (M, {\cal D}, \alpha)$.  The natural pairing
 \begin{equation}
 \ll (A_k, v_{k-1}), (c_k, \hat{X}, \sigma_{k-1}, \hat{Y}) \gg =  \langle  (\hat{X}^*(A_k), \hat{Y}^*(v_{k-1})),
 (c_k, \sigma_{k-1}) \rangle,
 \label{natrelabdow}\end{equation}
  defined through the pairing (\ref{naturalapeuska}), has the property
   \begin{equation}
    \ll \delta (A_k, v_{k-1}), (c_{k+1}, \hat{X}, \sigma_{k}, \hat{Y}) \gg
    =
     \ll (A_k, v_{k-1}), \d (c_{k+1}, \hat{X}, \sigma_{k}, \hat{Y}) \gg ,
   \label{stropcklwpw}\end{equation}
 which shows that
$\delta$ is adjoint to $\d$.

\end{document}